\pdfoutput=1
\documentclass[a4paper,10pt]{amsart}
\usepackage{amsmath}
\usepackage{amssymb}
\usepackage{amsxtra}
\usepackage{amsthm,amsfonts,latexsym}
\usepackage{epic}
\usepackage{eepic}
\usepackage{fancybox}
\usepackage{boxedminipage}
\usepackage{calc}
\usepackage{tabularx}
\usepackage{graphicx,color}
\usepackage{longtable,booktabs}
\usepackage{supertabular,longtable}
\usepackage{paralist}
\usepackage{amsmath}
\usepackage{amsfonts}
\usepackage{amssymb}
\usepackage{amsthm}
\usepackage[english]{babel}
\usepackage[latin1,ansinew]{inputenc}
\usepackage{a4wide}
\usepackage{pdfpages}

\usepackage{color}
\newtheorem{thm}{Theorem}%[section]

\newcommand{\bea}{\begin{eqnarray}}
\newcommand{\eea}{\end{eqnarray}}
\newcommand{\beann}{\begin{eqnarray*}}
\newcommand{\eeann}{\end{eqnarray*}}
\newcommand{\ba}{\begin{array}}
\newcommand{\ea}{\end{array}}

\newcommand{\beq}{\begin{equation}}
\newcommand{\eeq}{\end{equation}}
\newcommand{\be}{\begin{equation}}
\newcommand{\ee}{\end{equation}}

\newcommand{\eps}{\epsilon}

\accentedsymbol{\hcirc}{ {\overset{\scriptscriptstyle \circ }{ {\rm H} }}}

\newcommand{\R}{\mathbb{R}}

\newcommand{\divv}{ \nabla \!\! \cdot \!}
\newcommand{\bu}{{\bf u}}

\newcommand{\ovrho}{{\overline\rho}}
\newcommand{\unrho}{{\underline\rho}}
\newcommand{\ovc}{{\overline c}}
\newcommand{\unc}{{\underline c}}

\newtheorem{theorem}{Theorem}
\title{Diffuse planar phase boundaries in a two-phase fluid\\ with one very dense phase}
\author{Heinrich Freist\"uhler and Matthias Kotschote}
\date{July 11, 2013}
\begin{document}
\begin{abstract}
This note studies Navier-Stokes-Allen-Cahn models for 
compressible fluids that are mixtures of two incompressible phases 
whose density ratio $\eps=\rho_1/\rho_2$ is very small. Under a natural assumption
on the mixing energy, it shows the existence of 
diffuse planar phase boundaries for all $0\le\eps<\eps_0$. For $\eps=0$, one recovers  
the Navier-Stokes-Korteweg model and its well-known diffuse phase boundaries.
\end{abstract}
\maketitle
%\thispagestyle{empty}
%\numberwithin{equation}{section}
%%%%%%%%%%%%%%%%%%%%%%%%%%%%%%%%%%%%%%%%%%%%%%%%%%%%%%%%%%%%%%%%%%%%%%%%%%%%%
%%%%%%%%%%%%%%%%%%%%%%%%%%%%%%%%%%%%%%%%%%%%%%%%%%%%%%%%%%%%%%%%%%%%%%%%%%%%%
\parindent=0cm
We consider the Navier-Stokes-Allen-Cahn system of Blesgen \cite{B},
\begin{alignat}{3}
\partial_t \rho + \divv ( \rho \bu ) & = 0, & \quad 
\notag
\\
\partial_t (\rho \bu)  + \divv ( \rho \bu \otimes\bu + p{\bf I}) 
& = 
\divv \left(
\mu (\nabla\bu + (\nabla \bu)^T)
+
(\lambda \divv \bu) {\bf I} 
- \delta \rho \nabla c \otimes \nabla c 
\right), 
\label{NSAC}
\\
\partial_t (\rho c)  + \divv (\rho c\, \bu ) & =
\delta^{-1/2}(\rho q
+
\divv \left(\delta \rho \nabla c \right))
\notag
\end{alignat}
with 
$$
\rho^{-1}=\tau=G_p(p,c) \quad\hbox{and}\quad q=-G_c(p,c)
$$
and 
\begin{equation}\label{Gibbs}
G(p,c)
=
\hat G(p,c)
+W(c)
\end{equation}
with
\begin{align}
\label{LT}
\hat G(p,c)&=(c\tau_1+(1-c)\tau_2)p\\
&=(c+(1-c)\eps)\tau_*p\quad\hbox{with fixed }\tau_*>0.\label{FK} 
\end{align}
Note that \eqref{Gibbs}, \eqref{LT} is the Gibbs energy used by Lowengrub and Truskinovsky 
in their pioneering work on two-fluid mixtures \cite{LT}, in which they formulate and treat 
the corresponding Navier-Stokes-Cahn-Hillard system. Differently from their situation, ours 
allows phase transformations. \\

The results of this note are contained in the following two theorems.
\begin{thm}
Assume that $W\in C^2((0,1),(0,\infty))$ satisfies 
$$ %\begin{equation}
W(0+)>W(\overline{c})>W(\underline{c})>W(1-)
$$ %\end{equation}
and 
$$ %\begin{equation}
W'(c)\neq 0\quad \hbox{for all }c\in(0,1)\setminus\{\underline{c},\overline{c}\}
$$ %\end{equation}
with two concentrations
$$
0<\underline{c}<\overline{c}<1.
$$
If $\eps,|m|\ge 0$ are sufficiently small,
then there exists an undercompressive planar heteroclinic traveling wave solution of \eqref{NSAC} 
with  \eqref{Gibbs},\eqref{LT},\eqref{FK} that has mass transfer $m$.
\end{thm}

\begin{thm}
In the limiting case $\eps=0$, the system \eqref{NSAC} can 
be written as the Navier-Stokes-Korteweg system
\begin{alignat}{3}
\partial_t \rho + \divv ( \rho \bu ) & = 0, & \quad 
\notag\\
\partial_t (\rho \bu)  + \divv ( \rho \bu \otimes\bu )
- 
\divv (\mathcal{S}[\bu] + \mathcal{K}[\rho]) & = 0, & \quad
\label{NSK}
\end{alignat}
with
\begin{align}
\mathcal{S}[\bu] & = 2 \mu \mathcal{D}(\bu) + ( \lambda + \delta^{1/2}\tau_*^{-2}\rho^{-1} ) \divv \bu \, \mathcal{I} , 
\label{S}\\
\mathcal{K}[\rho] & = - \rho^2 \psi_\rho \, \mathcal{I} + \rho \divv ( \kappa(\rho) \nabla \rho ) \, \mathcal{I} 
- \kappa(\rho) \nabla \rho \otimes \nabla \rho,\label{K}
\\
\psi & = W((\tau_*\rho)^{-1}) + \frac{\kappa(\rho)}{2 \rho} |\nabla \rho|^2,
\quad \kappa(\rho) = \delta \tau_*^{-2}\rho^{-3}.
\label{psi}
\end{align}
\end{thm}
{\bf Proof.} For $\eps=0$, the specific volume 
$$ %\begin{equation}
\tau=c\tau_1+(1-c)\tau_2
$$ %\end{equation}
is 
$$
\tau=c\tau_*.
$$
Using this in $(1)_3$ yields 
$$
p
=
\tau_*^{-1}
\left[
-W'(c)  
+\rho^{-1}\divv \left(\delta \rho \nabla c \right)
\right]
-\delta^{1/2}\tau_*^{-2}\rho^{-1}\divv \bu. 
$$
Inserting this in $(1)_2$ results in \eqref{NSK} with \eqref{S}, \eqref{K}, \eqref{psi}. This proves
Theorem 2. For the limiting case $\eps=0$, Theorem 1 now follows from 
results of Slemrod \cite{S} and Benzoni-Gavage et al.\  \cite{BG}.
Considering the ODE system for traveling waves of \eqref{NSAC} as in \cite{FK}, the general result of 
Theorem 1 follows by perturbation.\\

We conclude by remarking that again along the lines of \cite{F3},
existence results for diffuse boundaries that are similar to those detailed here and in \cite{FK} 
hold for fluids consisting of two {\it compressible}\ immiscible phases.
In that case one considers for instance energies of the form
$$
U(\tau,c)=cU_1(\tau_1(\tau,c))+(1-c)U_2(\tau_2(\tau,c))
$$ 
composed from the two single-phase energies
$$
U_j:(0,\infty)\to(0,\infty) 
$$
with
$$
U_j(0+)=-U_j'(0+)=\infty,\quad U_j(\infty)=U_j'(\infty)=0
\quad
\hbox{and}
\quad
U_j''>0,\quad j=1,2,
$$
and $\tau_j(\tau,c), j=1,2,$ defined via
\begin{equation}\label{tau}
\tau=c\tau_1+(1-c)\tau_2
\end{equation}
and
\begin{equation}\label{p}
-U_1'(\tau_1)=-U_2'(\tau_2).
\end{equation}
While \eqref{tau} expresses immiscibility with zero excess volume,
the two values equated in \eqref{p} are identical with the common pressure 
$$
p=-U_\tau(\tau,c)
$$
of the two phases.

\end{document}